\documentclass[a4paper,10pt]{article}
\usepackage[T1]{fontenc}
\usepackage[latin1]{inputenc}
\usepackage[dvips]{graphicx}
\usepackage[dvips]{color}
\usepackage{amssymb}
\usepackage{amsfonts}
\usepackage{amsmath}
\usepackage[active]{srcltx}

\begin{document}


\title{Some comments concerning the blow-up of solutions of the exponential reaction-diffusion equation}
\author{A. Pulkkinen}
\date{}
\maketitle



\newtheorem{theo}{Theorem}
\newtheorem{pro}{Proposition}[section]
\newtheorem{lem}[pro]{Lemma}
\newtheorem{defi}[pro]{Definition}               
\newtheorem{coro}[pro]{Corollary}
\newtheorem{rema}[pro]{Remark}
\newtheorem{assu}[pro]{Assumption}

\newcommand{\ma}{\mathcal{A}}
\newcommand{\mb}{\mathcal{B}}
\newcommand{\mc}{\mathcal{C}}
\newcommand{\mf}{\mathcal{F}}
\newcommand{\mg}{\mathcal{G}}
\newcommand{\ml}{\mathcal{L}}
\newcommand{\mm}{\mathcal{M}}
\newcommand{\mn}{\mathcal{N}}
\newcommand{\mr}{\mathcal{R}}
\newcommand{\mt}{\mathcal{T}}
\newcommand{\mx}{\mathcal{X}}
\newcommand{\my}{\mathcal{Y}}
\newcommand{\mz}{\mathcal{Z}}                    
\newcommand{\mbc}{\mathbb{C}}
\newcommand{\mbn}{\mathbb{N}}
\newcommand{\mbz}{\mathbb{Z}}
\newcommand{\mbx}{\mathbb{X}}
\newcommand{\mbr}{\mathbb{R}}
\newcommand{\mbe}{\mathbb{E}}
\newcommand{\mbt}{\mathbb{T}}
\newcommand{\mbp}{\mathbb{P}}

\newcommand{\bdm}{\begin{displaymath}}
\newcommand{\edm}{\end{displaymath}}
\newcommand{\be}{\begin{equation}}
\newcommand{\ee}{\end{equation}}                   
\newcommand{\bmu}{\begin{multline*}}
\newcommand{\emu}{\end{multline*}}
\newcommand{\bea}{\begin{eqnarray*}}
\newcommand{\eea}{\end{eqnarray*}}

\newcommand{\comment}[1]{}
\newcommand{\eps}{\epsilon}
\newcommand{\veps}{\varepsilon}
\newcommand{\av}{\Arrowvert}
\newcommand{\ud}{\mathrm{d}}
\newcommand{\udm}{\mathrm{d} \mu}
\newcommand{\udmy}{\, \mathrm{d} \mu (y)}
\newcommand{\nelio}{\hspace{\stretch{1}}$\Box$}      
\newcommand{\wt}{\widetilde}
\newcommand{\wh}{\widehat}
\newcommand{\ol}{\overline}
\newcommand{\inti}{\int\displaylimits}
\newcommand{\lra}{\Leftrightarrow}
\newcommand{\llra}{\Longleftrightarrow}
\newcommand{\lto}{\longrightarrow}


\begin{abstract}
The aim of this paper is to refine some results concerning the blow-up of solutions of the exponential reaction-diffusion equation. We consider solutions that blow-up in finite time, but continue to exist as weak solutions beyond the blow-up time. The main result is that these solutions become regular immediately after the blow-up time. This result improves on that of Fila, Matano and Pol\'acik, who consider radially nonincreasing solutions, whereas no such assumption is needed here. Under certain additional assumptions we also obtain that the regularization is asymptotically selfsimilar.

Moreover, we consider the question of blow-up rate for radial solutions of the two-dimensional problem and prove that the blow-up is of type I, provided that the maximum of the solution is attained at the origin.
\end{abstract}

\section{Introduction}

In this paper we consider the solutions of the equation
\begin{equation} \label{eq1}
\left\{
\begin{array}{lll}
u_t = \Delta u + f(u),& \quad x \in \Omega, &t>0, \\
u = 0,& \quad  x \in \partial \Omega, & t>0, \\
u(x,0) = u_0(x) \ge 0, & \quad x \in \Omega, &
\end{array}
\right.
\end{equation}
where $\Omega = B(R) = \{x \in \mbr^N \, : \, |x| < R\}$ and the initial condition $u_0$ is nonnegative and in $C^1(\ol{\Omega})$. We will assume that
\be \label{f(u)}
f(u) = e^u,
\ee
unless otherwise mentioned. This is in contrast to the most frequently studied nonlinearity in the blow-up theory, namely $f(u) = u|u|^{p-1}$.

A solution blows up in finite time if there exists $T \in (0,\infty)$ such that $\av u(\cdot,t) \av_{\infty} < \infty$ for $t < T$ and 
$$
\limsup_{t \to T} \av u(\cdot,t) \av_{\infty} \to \infty,
$$
as $t \to T$. In this paper we wish to refine some earlier results concerning both blow-up and regularity of solutions of equation (\ref{eq1}).

There are several questions that have been raised related to the behavior of blow-up solutions. One fundamental question concerns the blow-up rate of a solution. To that end, blow-up is categorized in two classes with respect to its rate. Since it is natural to compare the equation (\ref{eq1}) to its ordinary differential equation counterpart $u' = f(u)$, it has become standard to say that blow-up is of type I if the blow-up rate of the solution of (\ref{eq1}) is the same as the blow up rate of the solution of $u' = f(u)$. When taking $f(u) = e^u$, type I blow-up corresponds to the inequalities 
$$
C_1 \le \log(T-t) + \av u(\cdot,t) \av_{\infty} \le C_2, \quad \text{ for every } t \in (0,T),
$$
for some constants $C_1$ and $C_2$. If blow-up is not of type I then it is said to be of type II.

A blow-up point is defined to be a point $x_0 \in \ol{\Omega}$ such that there exists a sequence $\{(x_n,t_n)\}_n \subset \Omega \times (0,T)$ satisfying $(x_n,t_n) \to (x_0,T)$ and $u(x_n,t_n) \to \infty$ as $n \to \infty$.

By the standard theory of parabolic regularity we know that a solution is classical for any $t \in (0,T)$. Even though the regularity in that sense is lost at the blow-up moment, we can still talk about weak solutions which may exists also after the blow-up time. To that end we give the following definition.
\begin{defi}\label{d:L1}
By an $L^1$-solution of (\ref{eq1}) on $[0,\mt]$ we
mean a function $u\in C([0,\mt];$  $L^1(\Omega))$ such that $f(u)\in
L^1(Q_\mt)$, $Q_\mt:=\Omega\times(0,\mt)$ and the equality
$$\int_\Omega[u\Psi]^{t_2}_{t_1}\,dx-\int^{t_2}_{t_1} \int_\Omega u\Psi_t\,
dx\,dt=\int^{t_2}_{t_1} \int_\Omega (u\Delta\Psi+f(u)\Psi)\,dx\,dt$$
holds for any $0\le t_1<t_2\le \mt$ and $\Psi\in C^2(\bar Q_\mt)$,
$\Psi=0$ on $\partial\Omega\times [0,\mt]$.
\end{defi}
If a solution does not exist as a weak solution after the blow-up time, then  the blow-up is said to be complete.

In this paper we will discuss both the blow-up rate of certain solutions of (\ref{eq1})-(\ref{f(u)}) and the regularity of certain weak solutions after the blow-up time. Let us first describe our approach in proving regularity for weak solutions.

To consider the regularity of weak solutions after the blow-up time, we use the available information about the blow-up profile at the blow-up moment. The following Theorem describes the asymptotics of radially symmetric weak solutions that blow-up with type I rate with respect to the similarity variables as the blow-up moment is approached, see \cite{FP}, \cite{P}. The result is analogous to the results in \cite{M}, where solutions with nonlinearity $f(u) = u^p$ and $p\in (p_S,p_L)$ were discussed. Here the exponent
$$
p_{S} = \left\{
\begin{array}{ll}
\infty,& \text{ if } N \le 2, \\
\frac{N+2}{N-2}, & \text{ if } N > 2,
\end{array} \right.
$$
and the exponent $p_L = \frac{N-4}{N-10}$ for $N>10$ and $p_L = \infty$ for $N \le 10$. See also \cite{B}, \cite{BB}, \cite{GK}, \cite{HV}, \cite{Ve1} and \cite{Ve2} for further results concerning similarity variables and blow-up profiles of solutions.

\begin{theo} \label{theo:FP}
Let $u$ be a radially symmetric $L^1$-solution of (\ref{eq1})-(\ref{f(u)}) on $[0,\mt]$ that blows up with type I rate at $(x,t) = (0,T)$, where $T < \mt$. Then 
\begin{equation} \label{eq:ssconv1}
\lim_{t \uparrow T} \left( \log(T-t) + u(y\sqrt{T-t}, t)\right) = \varphi(y),
\end{equation}
uniformly for $y$ in compact sets of $\mbr^N$, where $\varphi$ satisfies
\begin{equation} \label{eq:stationary}
\left\{
\begin{array}{ll}
\Delta \varphi - \frac{y}{2} \nabla \varphi + e^\varphi -1 =0,& \qquad |y| >0, \\
\varphi(0) = \alpha, \; \nabla \varphi(0) = 0,
\end{array}
\right.
\end{equation}
and
\begin{equation} \label{stationaryAsympt}
\lim_{|y| \to \infty} \big( \varphi(y) +2\log |y| \big) = C_\alpha,
\end{equation}
for some $\alpha > 0$ and $C_{\alpha} \in \mbr$.
\end{theo}

\begin{rema} \label{rema_BU_origin}
By Corollary \ref{BU_origin} below the origin is the only blow-up point provided that $u$ is a radially symmetric minimal $L^1$-solution on $[0,\mt]$.
\end{rema}

The idea of the proof is to assume that the convergence (\ref{eq:ssconv1}) holds with $\varphi \equiv 0$ and then prove that either
\be \label{loglog_profile}
u(x,T) + 2\log|x| - \log|\log|x|| \to C, \quad \text{as } |x| \to 0,
\ee
or
$$
u(x,T) - m\log|x| \to C, \quad \text{as } |x| \to 0,
$$
for some constant $C$, which implies complete blow-up by the results in \cite{Va}. So, as a byproduct of the above Theorem \ref{theo:FP}, we obtained the final time blow-up profiles of solutions that blow-up with type I rate and have a constant selfsimilar blow-up profile, i.e., the convergence (\ref{eq:ssconv1}) holds with $\varphi \equiv 0$.

The following result is proved in a recent paper \cite{P}. There the final time blow-up profiles were found for solutions such as in the previous Theorem \ref{theo:FP}, that is, for solutions having a nonconstant selfsimilar blow-up profile. In \cite{MM2} similar results were proved for the nonlinearity $f(u)=u|u|^{p-1}$ and $p>p_S$.

\begin{theo} \label{theorem2}
Assume that $u$ is a solution of (\ref{eq1})-(\ref{f(u)}) that blows up at $(x,t) = (0,T)$ with type I rate and verifies (\ref{eq:ssconv1})-(\ref{stationaryAsympt}). Then the final time blow up profile of $u$ is given by
\be \label{eq:u_profile1}
|u(x,T) +2\log|x| - C_{\alpha}| \to 0, \quad \text{ as } |x| \to 0,
\ee
where $C_{\alpha}$ is the constant from (\ref{stationaryAsympt}).
\end{theo}

The above Theorem allows us to use a result in \cite{Va} to attack the question of regularity after the blow-up. However, we have to work with weak solutions that are limits of classical solutions and so let us give the following definitions.
\begin{defi}
By a limit $L^1$-solution of (\ref{eq1}) on $[0,\mt)$ we mean a function $u$ that satisfies the following. It can be approximated by a sequence $\{u_n\}_{n=1}^{\infty}$ of functions such that $u_n$ verifies (\ref{eq1}) on $[0,\mt)$ in the classical sense with initial data $u_{0,n} \in C(\ol{\Omega})$ and
$$
u_{0,n} \to u_0, \quad \text{in } C(\ol{\Omega}),
$$
and
\begin{align*}
&u_n \to u, \quad \text{in } L^1(\Omega) \text{ for every } t\in [0,\mt), \\
&f(u_n)\to f(u), \quad \text{in } L^1(\Omega \times (0,t)) \text{ for every } t\in [0,\mt).\\
\end{align*}

A limit $L^1$-solution is said to be a minimal $L^1$-solution if the approximating initial data $\{u_{0,n}\}_n$ verify
$$
0 \le u_{0,1}(x) \le u_{0,2}(x) \le u_{0,3}(x) \le \ldots ,
$$
for every $x \in \Omega$.
\end{defi}

It follows from the definitions of a limit $L^1$-solution and a minimal $L^1$-solution that these are also $L^1$-solutions.

The main result of this paper is the following.

\begin{theo} \label{theorem3}
Let $ \in [3,9]$ and assume that $u$ is a radially symmetric minimal $L^1$-solution of (\ref{eq1})-(\ref{f(u)}) on $(0,\mt)$ that blows up at $t = T < \mt$ with type I rate and that the convergence (\ref{eq:ssconv1})-(\ref{stationaryAsympt}) takes place. 

Then $u$ is regular, i.e. $u(\cdot,t) \in C^{\infty}(\Omega)$, on some time interval $(T,T+\eps)$.

Furthermore, it holds that $C_{\alpha} \le c^{\#}$, where $c^{\#}$ depends only on the dimension of the space, and if $C_{\alpha} < c^{\#}$ or $C_{\alpha}=c^{\#}$ but $u(x,T) < -2\log|x| + c^{\#}$ near the origin, then the regularization is selfsimilar. This means that the rate of regularization is given by
\be \label{eq:reg_rate}
\log(t-T) + \av u(\cdot,t) \av_{\infty} \le C,
\ee
for some contant $C>0$ and $t \in (T,T+\eps)$. Moreover, it holds 
\be \label{eq:reg_conv}
\lim_{t \downarrow T} \left(\log(t-T) + u(\sqrt{t-T}y,t) \right) = \psi(y),
\ee
uniformly for $y$ in compact sets, where
\be \label{eq:fw_selfs}
\left\{
\begin{array}{ll}
\Delta \psi + \frac{y}{2} \nabla \psi + e^{\psi} + 1 =0, & \qquad |y| > 0, \\
\psi(0) = \beta, \nabla \psi(0) = 0. 
\end{array}
\right.
\ee
\end{theo}

\begin{rema} \label{remark}
By Theorem \ref{theo:FP} and Remark \ref{rema_BU_origin}, the convergence (\ref{eq:ssconv1})-(\ref{stationaryAsympt}) takes place provided that $u$ is a radial and minimal $L^1$-solution of (\ref{eq1})-(\ref{f(u)}) on $(0,\mt)$ that blows up at $t = T < \mt$ with type I rate.

Moreover, the blow-up is of type I if $N \in [3,9]$ and $u(0,t) = \max_{x\in \Omega}u(x,t)$ by Theorem 1.1 in \cite{FP}.
\end{rema}

Immediate regularization after the blow-up has been proved in \cite{FMP} for $f(u) = e^u$ under the extra assumption that $u$ is radially nonincreasing. Here we are able to slightly improve the result by not assuming that $u$ is radially nonincreasing. We also obtain the rate of regularization, except for the one special case.

In the case $C_{\alpha} < c^{\#}$ or $C_{\alpha} = c^{\#}$ but $u(x,T) < -2\log|x| + c^{\#}$ for $x$ close to the origin, we can use techniques from \cite{Va} to obtain the result of the above Theorem. In this case the approach is also independent of the assumption on radial symmetry. If $u(x,T)$ may attain larger values than $-2\log|x| + c^{\#}$ for $x$ close to zero, we have to refine the approach of \cite{FMP} to obtain the claim.

Questions related to blow-up profiles of solutions and regularization after the blow-up are discussed in the paper \cite{MM2} for $f(u) = u|u|^{p-1}$. Immediate regularization for any limit $L^1$-solution is proved there, assuming only supercriticality and type I blow-up. Our incapability of dealing with nonminimal continuations derives from the techniques used to prove Theorem \ref{theorem2} as we do not obtain apriori bounds for the solutions of (\ref{eq1}).

Nonuniqueness of $L^1$-continuations of $u$ for $f(u) = u^p$ was proved in \cite{FM}.

In the last section of this paper we will consider the blow-up rate of solutions of (\ref{eq1})-(\ref{f(u)}). Type I blow-up is a frequent phenomenon when considering equation (\ref{eq1}) and in fact when $f(u) = u|u|^{p-1}$ every blow-up solution exhibits type I blow-up in the subcritical case, that is, for $p \in (1,p_S)$.

Blow-up is also of type I in the supercritical case, if $u$ is radially symmetric and $p \in (p_S,p_{JL})$, where
$$
p_{JL} =
\left\{
\begin{array}{ll} \infty, & n \le 10 \\
\frac{4}{n-4-2\sqrt{n-1}},& n > 10.
\end{array}
\right.
$$
These results are due to the classical papers \cite{GK} and \cite{GMS} in the so-called subcritical range and \cite{MM1} when $p$ is supercritical.

For the exponential nonlinearity the results remain incomplete. Blow-up is known to be of type I in the subcritical case if either $N =1$ or $N=2$, the solution is radially symmetric and radially decreasing, \cite{HV3}, \cite{FHV}. Moreover, type I blow-up is known always to occur for $N \in (2,10)$ when the solution is radially symmetric and $u$ attains its maximum at the origin, see \cite{FP}.

In this paper we prove the following.
\begin{theo} \label{theo:n=2}
Let $u$ be a radially symmetric solution of (\ref{eq1})-(\ref{f(u)}) with $N = 2$ and $\Omega = B(R)$ or $\Omega = \mbr^2$, and assume that $u(0,t) = \max_x u(x,t)$. If $u$ blows up at $t = T<\infty$, then the blow-up is of type I.
\end{theo}

This result improves somewhat the result in \cite{FMP}, where it is assumed that $u$ is radially decreasing. However, we have to impose the assumption that $u$ attains its maximum at the origin. This assumption arises from the form of the intrinsic rescaling associated to the exponential equation which does not preserve positivity, in contrast to the power case.

The above Theorem \ref{theo:n=2} could also be proved by using the technique from \cite{GP} which is also used in \cite{FHV}. However, that approach is vitally dependent on the assumption that $u$ attains its maximum at the origin. Our approach depends on the assumption through Lemma \ref{lemma_conv}. By finding suitable estimates for $u$, Lemma \ref{lemma_conv} could presumably be proven also without assuming $u(0,t) = \max_x u(x,t)$, but we are unable to provide a proof here.

Even though type I blow-up is very common, also type II blow-up is known to take place. It is proved in \cite{HV2} and \cite{Mi1} that for $f(u)= u|u|^{p-1}$ there exist solutions that blow-up with type II rate when $p > p_{JL}$.

In the next section, we will first consider the proof of Theorem \ref{theorem3} when $u(x,T) < -2\log|x| + c_{\#}$ near the origin, and then, in Section \ref{section_apriori} we will finalize the proof by achieving some apriori regularity for the solution after the blow-up time.

The final section of this paper is devoted to the proof of Theorem \ref{theo:n=2}.

\section{Regularization after blow-up} \label{section_regul}

In this section we prove Theorem \ref{theorem3}. First we prove one part of the Theorem as a consequence of Theorem \ref{theorem2} and of the following result proved in \cite{Va}.

\begin{pro} \label{pro:Va}
Consider the equation (\ref{eq1})-(\ref{f(u)}) with $\Omega = B(R)$ and $N \in [3,9]$. There exists a constant $c^{\#}>\log(2(N-2))$ such that if $u$ is a minimal $L^1$-solution of the problem that blows up at $t=T$ and 
$$
u(x,T) < -2\log|x| + c^{\#}
$$
in a neighborhood of $x = 0$, then $u$ is regular immediately after the blow-up and 
$$
\log(t-T) + \av u(\cdot,t) \av_{\infty} < C,
$$
for $t>T$ close enough to $T$.

If on the other hand 
$$
u(x,T) \ge -2\log|x| + c,
$$
with some $c > c^{\#}$ in a neighborhood of $x = 0$, then the blow-up is complete.
\end{pro}

Assume now that $u$ is as in Theorem \ref{theorem3}, which implies by Theorem \ref{theorem2} that
$$
\lim_{x \to 0} (u(x,T) + 2\log|x|) = C_{\alpha}.
$$

If $C_{\alpha} > c^{\#}$, then $u$ would not be an $L^1$-solution by Proposition \ref{pro:Va} and so $C_{\alpha} \le c^{\#}$ as stated in Theorem \ref{theorem3}.

If $C_{\alpha} < c^{\#}$ or $C_{\alpha} = c^{\#}$ and $u(x,T) < 2\log|x| + c^{\#}$ near the origin, then Proposition \ref{pro:Va} gives the first part of the claim.

Let us now prove the convergence to a forward selfsimilar solution as $t$ approaches $T$ from above. We prove this by comparing forward selfsimilar solutions with the solution $u$.

It is also proved in \cite{Va} that if we take initial data $u_0(x) = -2\log|x| + c$ with $c \le c^{\#}$, then the unique so-called proper solution (see \cite{Va} for further details) of (\ref{eq1})-(\ref{f(u)}) is the forward selfsimilar solution given by
$$
u_c(x,t) = -\log(t) + \psi\left(\frac{x}{\sqrt{t}}\right),
$$
where $\psi$ satisfies (\ref{eq:fw_selfs}) with $\lim_{|y|\to \infty} (\psi(y) + 2\log|y|)= c$.

Let $u_n(x,t)$ be the approximating sequence for the minimal $L^1$-solution. Then, we have that $u_n(x,T) \le u(x,T) < -2\log|x| + C_{\alpha}$ for $|x| \le 2\eps_1$ for some $\eps_1 > 0$ independent of $n$. Therefore we can prove, as was done in the proof of Proposition \ref{pro:Va} in \cite{Va}, that $u(x,t) < u_{C_{\alpha}}(x,t-T)$ for $|x| \le \eps_1$ and $t \in (T,T+\delta_1)$ for some $\delta_1 > 0$. 

Let $u_{0,n,c'}$ be the approximating sequence related to the proper solution with initial data $u_{0,c'} = -2\log|x| + c'$ where $c' < C_{\alpha}$. Then $u_{0,n,c'}(x)  < -2\log|x| + c' < u(x,T)$ for $|x| < 2\eps$ for some $\eps > 0$ independent of $n$, and so we can prove that $u_{c'}(x,t-T) < u(x,t)$ for $|x| < \eps$ and $t \in (T,T+\delta)$ for some $\delta > 0$, again using the same method as in \cite{Va}. We have thus obtained that for every $\sigma > 0$ there exists $\eps_2(\sigma),\delta_2(\sigma) > 0$ such that
\be \label{forward_esti}
u_{C_{\alpha} - \sigma}(x,t-T) < u(x,t),
\ee
for every $t \in (T,T+\delta_2(\sigma))$ and $|x| < \eps_2(\sigma)$.

Let $\theta$, $M > 0$ and take $\sigma>0$ small enough so that
$$
|\psi_{C_{\alpha}-\sigma}(y) - \psi_{C_{\alpha}}(y)| < \theta,
$$
for every $|y| < M$. This gives that by taking $\delta_3 \in (0,\min\{\delta_1,\delta_2(\sigma)\})$ small enough such that $\sqrt{\delta_3} M < \min\{\eps_1,\eps_2(\sigma)\}$, we have
\begin{multline*}
|\log(t-T) + u(\sqrt{t-T}y,t) - \psi_{C_{\alpha}}(y)|
\\
= |u(\sqrt{t-T}y,t) - u_{C_{\alpha}}(\sqrt{t-T}y,t-T)|
\\
\le |u_{C_{\alpha}-\sigma}(\sqrt{t-T}y,t-T) - u_{C_{\alpha}}(\sqrt{t-T}y,t-T)|
\\
= |\psi_{C_{\alpha}-\sigma}(y) - \psi_{C_{\alpha}}(y)| < \theta,
\end{multline*}
for every $|y| \le M$ and $t \in (T,T+\delta_3)$.

Therefore,
$$
\lim_{t \downarrow T} |\log(t-T) + u(\sqrt{t-T}y,t) - \psi_{C_{\alpha}}(y)| = 0,
$$
and we have proved the selfsimilarity of the regularization, i.e., that (\ref{eq:reg_rate})-(\ref{eq:fw_selfs}) hold in the case $C_{\alpha} \le c^{\#}$ and $u(x,T) < -2\log|x| + c^{\#}$ near the origin.

The part concerning $C_{\alpha} = c^{\#}$ and $u(x,T)$ attaining also larger values than the function $-2\log|x| + c^{\#}$ near the origin will be discussed in the following.

\subsection{Apriori regularity of the continuation} \label{section_apriori}

Here we want to demonstrate that the minimal limit $L^1$-continuation has some apriori regularity also after the blow-up time and conclude the proof of Theorem \ref{theorem3}. We will proceed along the lines of \cite{FMP}, but make some modifications in order to be able to handle radially nondecreasing solutions as well.

First we will show, in Lemma \ref{lemma_C_alpha} and Lemma \ref{lemma_C_C}, that the solution blows up only at the origin and that it is continuous in time, and locally Hölder continuous outside the origin. Then we will use an intersection comparison method in Lemma \ref{pro:loglog} to prove an upper bound for the solution, which is valid for $t$ strictly larger than the blow-up moment. In Lemma \ref{lem:BUtimes} we prove that the solution can blow-up at most finitely many times. Then, we proceed by showing in Lemma \ref{lem:decreasing} that if the solution does not become regular immediately after the blow-up, then it is radially nonincreasing close to the blow-up point $x=0$. These results allow us to use the same techniques as in \cite{FMP} in Lemma \ref{lem:C2reg}, Lemma \ref{Lemma_log} and Proposition \ref{pro:reg} to conclude the proof.

Since we will work with a radially symmetric function $u$, it will be useful to use radial notation. So, define
\be \label{Udef}
U(|x|,t) = u(x,t).
\ee
We will frequently use the zero number diminishing property for the solutions of one dimensional parabolic equations, so let us state it in the next Proposition, see Theorem 52.28 in \cite{QS}.
\begin{pro} \label{pro:zeron}
Let the zero number of a function $v \in C((0,R))$ be defined through
$$
\begin{array}{ll}
z_{[0,R]}(v) = \sup\left\{k \in \mbn \, : \right. & \text{there exist } 0 < x_1 < \ldots < x_k < R \\
 & \left. \text{such that } v(x_i)v(x_{i+1}) < 0, \text{ for } 0 \le i \le k-1\right\}.
\end{array}
$$
Assume that $V$ is a nontrivial classical solution of the equation
\be \label{V_eq}
V_t = V_{rr} + \frac{N-1}{r}V_r + QV, \quad \text{ for } r\in (0,R), \, t \in (t_1,t_2),
\ee
for $N \ge 1$, with $V_r(0,t)=0$ for $t \in (t_1,t_2)$ and either $V(R,t)=0$ or $V(R,t) \ne 0$ for every $t \in (t_1,t_2)$, where $Q \in L^{\infty}([0,R]\times (t_1,t_2))$.

Then 
\begin{enumerate}
\item $z_{[0,R]}(V(\cdot,t))< \infty$, for every $t \in (t_1,t_2)$. \\
\item $z_{[0,R]}(V(\cdot,t))$ is nonincreasing in time.\\
\item If $V(r_0,t_0) = V_r(r_0,t_0)=0$ for some $(r_0,t_0) \in [0,R] \times (t_1,t_2)$, then
$$
z_{[0,R]}(V(\cdot,t')) < z_{[0,R]}(V(\cdot,t'')),
$$
for any $t_1<t'<t_0<t''<t_2$.
\end{enumerate}

The same conclusion holds also if $N=1$ and $V$ is considered on an interval $[R_1,R_2]$ and $V(R_i,t) = 0$  or $V(R_i,t) \ne 0$ for every $t \in (t_1,t_2)$.
\end{pro}

The above property generalizes in the context of limit $L^1$-solutions as follows, see Lemma 4.1 in \cite{FMP}.
\begin{rema}
Results analogous to those of Proposition \ref{pro:zeron} can be formulated also if
$$
V(|x|,t) =u(|x|,t)-\varphi^*(|x|),
$$
where $u$ is a limit $L^1$-solution of the equation (\ref{eq1})-(\ref{f(u)}) on $[0,\mt)$ and
$$
\varphi^*(x)=-2\log|x|+ \log(2(N-2))
$$
is the singular solution of the equation (\ref{eq1})-(\ref{f(u)}). In this case the conclusion is the following.

Assume that there exist $0 < t_1 < s_1 < s_2 < \ldots s_k < t_2 < \mt$ and $\{r_i\}_{i=1}^k \subset (0,R]$ such that $V_r(r_i,s_i) = V(r_i,s_i)=0$. Then
$$
z_{[0,R]}(V(\cdot,t)) \le z_{[0,R]}(V(\cdot,t_1)) - k,
$$
for every $t \in (s_k,t_2)$. If no such $s_i$ can be found, then $k=0$. 
\end{rema}

The next Lemma states that the radial symmetry and minimality of the continuation actually implies some regularity for the solution $u$ of (\ref{eq1})-(\ref{f(u)}) after the blow-up time.
\begin{lem} \label{lemma_C_alpha}
Let $u$ be a radially symmetric minimal limit $L^1$-solution of (\ref{eq1})-(\ref{f(u)}) on $[0,\mt)$ that blows up at $t = T < \mt$. Then there exists $\alpha > 0$ such that $U(\cdot,t) \in C^{\alpha}_{\text{loc}}((0,R))$ for every $t \in (0,\mt)$, where $U$ is as in (\ref{Udef}).
\end{lem}
\emph{Proof.} Let $\{u_n\}_n$ be the approximating sequence related to the minimal solution $u$. Then $u_n(x,t) \le u(x,t)$ for every $t \in [0,\mt)$ and $x \in B(R)$. 

Standard regularity theory for the heat semigroup gives us the estimates
$$
\av e^{t \Delta}\phi\av_{L^q} \le (4 \pi t)^{-N\alpha/2} \av \phi\av_{L^p},
$$
for any $\phi \in L^p(B(R))$, where $1 \le p < q \le\infty$ and $\alpha = \frac 1p - \frac 1q$, and 
$$
\av e^{t \Delta }\phi \av_{X_{\beta}} \le C_{\beta} t^{-\beta} \av \phi \av_{L^q},
$$
for $\phi \in L^q(B(R))$, where $X_{\beta}=W^{2 \beta, q}_{0}(B(R))$ with $\beta \in (0,1)$ and $q \in (1,\infty)$.

By using these estimates with $\beta = \frac{1}{2}$, $q \in (1,\frac{N}{N-1})$, $p=1$ and $\alpha= 1 - \frac 1q \in (0, \frac 1N)$, we have the inequality
\be \label{Lap_X_b}
\av e^{t \Delta}\phi \av_{X_{1/2}} \le C \left( \frac t2 \right)^{-1/2} \av e^{t \Delta /2} \phi \av_{L^q} \le C \left( \frac t2 \right)^{-1/2} (2 \pi t)^{-N\alpha/2} \av \phi \av_{L^1}.
\ee

Because $u_n$ is a classical solution, it satisfies the variation of constants formula for every $t < \mt$. Therefore, we can use the above estimate (\ref{Lap_X_b}) to give us that
\begin{multline*}
\av u_n(\cdot,t+\delta) \av_{X_{1/2}} \le \av e^{\delta \Delta} u_n(\cdot,t)\av_{X_{1/2}} + \int_{t}^{t+\delta} \av e^{(t+\delta - \tau)\Delta} e^{u_n(\cdot,\tau)} \av_{X_{1/2}} \ud \tau 
\\
\le C(\alpha) \delta^{-1/2 - N\alpha/2} \av u_n(\cdot,t) \av_{L^1} + \int_{t}^{t+\delta} C(\alpha) (t+\delta-\tau)^{-1/2-N\alpha/2} \av e^{u_n(\cdot,\tau)} \av_{L^1} \ud \tau,
\end{multline*}
where $-1/2 - N\alpha/2 > -1/2 -1/2 = -1$. Therefore, for $0<\delta<t+\delta\le T_1<\mt$, it holds
\begin{multline*}
\av u_n(\cdot,t+\delta) \av_{X_{1/2}} 
\\
\le C(\alpha) \delta^{-(1+ N\alpha)/2} \sup_{\tau \in (0,T_1)} \av u(\cdot,\tau) \av_{L^1} + C(\alpha) \frac{\delta^{-(1 + N\alpha)/2 + 1}}{-(1 + N\alpha)/2 + 1} \sup_{\tau \in (0,T_1)} \av e^{u(\cdot,\tau)} \av_{L^1}
\\
< C(\alpha,\delta,T_1,u),
\end{multline*}
where we used the fact that $u_n(x,t) \le u(x,t)$ for $(x,t) \in B(R)\times [0,\mt)$.

We have thus obtained that $\{u_n(\cdot,t)\}_n$ is contained in a compact subset of the space $W^{1,q}_0(B(R))$ for every $t \in (\delta,T_1)$. Therefore, we can take a subsequence converging to a function $\wh{u}(\cdot,t) \in W_0^{1,q}(B(R))$ and, by the $L^1$-convergence of $u_n(\cdot,t)$ to $u(\cdot,t)$, it has to hold that $\wh{u}(\cdot,t) = u(\cdot,t)$ and that $u_n(\cdot,t) \to u(\cdot,t)$ in $W_0^{1,q}(B(R))$ along the original sequence for every $t \in (\delta,T_1)$.

Now we can use the radial symmetry of $u$ and Sobolev embedding in one dimension. To be precise, we know that the weak derivative $U_r = \nabla u \cdot \frac{x}{|x|}$ exists and
\begin{multline*}
\av U(\cdot,t) \av^{q}_{W^{1,q}([R_1,R_2])} \le C' \int_{R_1}^{R_2} |U(r,t)|^{q} + |U_r(r,t)|^{q} \ud r
\\
\le C' R_1^{1-N} \int_0^R |U(r,t)|^{q} + |U_r(r,t)|^{q} r^{N-1} \ud r \le C R_1^{1-N} \av u(\cdot,t) \av_{W^{1,q}(B(R))}^{q},
\end{multline*}
for $t \in (\delta,T_1)$. Hence $U(\cdot,t) \in W^{1,q}([R_1,R_2])$ for every $t \in (\delta,T_1)$ and $0<R_1<R_2 \le R$. Because of the one dimensional Sobolev embedding
$$
W^{1,q}([R_1,R_2]) \subset C^{\alpha}([R_1,R_2])
$$
we have
\begin{multline} \label{C_alpha_estimate}
\av U(\cdot,t) \av_{C^{\alpha}([R_1,R_2])} \le C \av U(\cdot,t) \av_{W^{1,q}([R_1,R_2])}
\\
\le C' R_1^{1-N} \av  u(\cdot,t) \av_{W^{1,q}(B(R))} \le C' R_1^{1-N} C'(\alpha,\delta,T_1,u),
\end{multline}
for every $t \in (\delta,T_1)$, which is the claim. \nelio

\begin{coro} \label{BU_origin}
Let $u$ be a radially symmetric minimal limit $L^1$-solution of (\ref{eq1})-(\ref{f(u)}) on $[0,\mt)$ that blows up at $t=T<\mt$. Then $x=0$ is the only blow-up point.
\end{coro}

As consequence of Lemma \ref{lemma_C_alpha} and standard parabolic estimates we have the following.
\begin{lem} \label{lemma_C_C}
For $u$ as in Lemma \ref{lemma_C_alpha}, it holds that
$$
U \in C^1((0,\mt),C([R_1,R_2])),
$$
for every $0<R_1<R_2 < R$.
\end{lem}

The idea of the next Lemma is to verify that the solution cannot oscillate between very large values and much smaller values near the origin. Proving this then gives us an upper bound for the solution. 
\begin{lem} \label{pro:loglog}
Let $u$ be a radially symmetric minimal limit $L^1$-solution of (\ref{eq1})-(\ref{f(u)}) on $(0,\mt)$ that blows up at $t=T<\mt$. Then for any $\eps >0$ small enough, there exist $T_1, T_2 \in (T,T+\eps)$ and a constant $C > 0$ such that
$$
u(x,t) \le -2\log|x| + \log|\log|x|| + C,
$$
for every $t \in (T_1,T_2)$ and $|x| \in (0,R)$.
\end{lem}

\emph{Proof.} By the above Lemma \ref{lemma_C_alpha} we know that $u$ is continuous outside the origin for every $t \in (0,\mt)$.

We will first show that there exist constants $C>0$ and $\eps>0$ such that $U(r,t)$, for any $t \in (T,T+\eps)$, can have at most finitely many intersections with the function
$$
\omega_C(r)=-2\log(r) + \log|\log(r)| + C,
$$
for $r$ close to zero.

Take a solution $u^*$ of (\ref{eq1})-(\ref{f(u)}) that blows up only at $(x,t) = (0,T^*)$, and satisfies
$$
\lim_{|x|\to 0} \left(u^*(x,T^*) - \omega_{C^*}(|x|)\right) = 0,
$$
for some $C^* \in \mbr$. The existence of such solutions is proved in \cite{B}. Define $C = C^* + 1$ and assume by contradiction that for any $\eps>0$ there exists $t_0 \in (T,T+\eps)$ such that $U(\cdot,t_0)$ crosses the function $\omega_C$ infinitely many times near the origin, which implies that $U(r_i,t_0) > \omega_C(r_i)$ for some sequence $\{r_i\}_i$ tending to zero as $i \to \infty$. Now let $V^*(|x|,t) = u^*(x,t + T^* - t_0)$, which implies that $V^*$ blows up at $t = t_0$. By taking $\eps$ small enough, we can assume that $t_0 - T^* \in (0,T)$ and so there exists a finite $M$ such that
$$
z_{[0,R]}(U(\cdot,t_0-T^*)-V^*(\cdot,t_0-T^*)) = M.
$$

Then take a sequence $\{\rho_i\}_i$ tending to zero, such that $U(\rho_i,t_0) < -2\log(\rho_i) + c^{\#}+1$. The existence of such a sequence is guaranteed by Proposition \ref{pro:Va}, and by the fact that $u$ does not blow-up completely at $t = t_0$. Because $U(\cdot,t_0)$ is continuous outside the origin and since $U_n(\cdot,t_0) \to U(\cdot,t_0)$ in $L^1(B(R))$, where $U_n(|x|,t_0) = u_n(x,t_0)$ is the approximating sequence for $u$, for every $n$, there exist sequences $\{r_i^n\}_i$ and $\{\rho_i^n\}_i$ tending to zero satisfying $\rho_i^n \in (r_{i+1}^n,r_i^n)$ and 
$$
U_n(r_i^n,t_0) > \omega_C(r_i^n), \quad \text{ and } \quad U_n(\rho_j^n,t_0) < -2\log(\rho_j^n) + c^{\#}+1,
$$
for $i \in \{1,\ldots,l(n)\}$ and $j \in \{1,\ldots,k(n)\}$, where $l(n)$ and $k(n)$ tend to infinity as $n \to \infty$. Consequently, there exist $l_0$ and $k_0$ such that
\begin{multline*}
U_n(r_i^n,t_0) - V^*(r_i^n,t_0) = U_n(r_i^n,t_0) - \omega_{C^*}(r_i^n) + \omega_{C^*}(r_i^n) - V^*(r_i^n,t_0)
\\
> 1 + \omega_{C^*}(r_i^n) - V^*(r_i^n,t_0) > \frac 12,
\end{multline*}
for $i \in \{l_0, \ldots, l(n)\}$ and similarly
\begin{multline*}
U_n(\rho_j^n,t_0) - V^*(\rho_j^n,t_0) = U_n(\rho_j^n,t_0) - \omega_{C^*}(\rho_j^n) + \omega_{C^*}(\rho_j^n) - V^*(\rho_j^n,t_0)
\\
< -\log|\log(\rho_j^n)| + c^{\#}+1 - C^* + \frac{1}{2} < 0,
\end{multline*}
for $j \in \{k_0, \ldots,k(n)\}$. Therefore $U_n(\cdot,t_0)$ intersects the function $V^*(\cdot,t_0)$ at least $\min\{l(n)-l_0,k(n)-k_0\}$ times. Because, by the results in \cite{B}, we may assume that $V^*(\cdot,t)$ blows up only at the origin, we obtain, by the parabolic estimates, that $V^* \in C([R_1,R_2]\times[t_0-T^*,t_0])$ for every $R_1, R_2 \in (0,R)$. Since $U_n$ is a classical solution, one has that
$$
z_{[0,R]}(U_n(\cdot,t_0-\eps_n)-V^*(\cdot,t_0-\eps_n)) \ge \min\{l(n)-l_0,k(n)-k_0\},
$$
for some $\eps_n > 0$ small enough.

Moreover, since $u_n(\cdot,t_0-T^*) \to u(\cdot,t_0-T^*)$ in $C^2(B(R))$, we know that
$$
z_{[0,R]}(U_n(\cdot,t_0-T^*)-V^*(\cdot,t_0-T^*)) = M,
$$
whenever $n$ is large enough. By taking $n$ large enough such that $\min\{l(n)-l_0,k(n)-k_0\} > M$, we have a contradiction with the zero number property. Therefore, we have proved what we claimed, that is, $U(r,t)$ intersects the function $\omega_C(r)$ at most finitely many times for $r$ close to $0$ for every $t \in (T,T+\eps)$ and $\eps > 0$ small enough.

Let us prove the original claim. Since $U(\cdot,t)$ has only finitely many intersections with the function $\omega_C$ near the origin, we define $r = \delta(t)>0$ to be the first intersection of $U(r,t)$ and $\omega_C(r)$, for every $t \in (T,T+\eps)$. Then $U$ lies either above or below $\omega_C$ for $r < \delta(t)$ and by Proposition \ref{pro:Va} it has to hold that $U(r,t) \le \omega_C(r)$ for every $r \le \delta(t)$ and $t \in (T,T+\eps)$. 

We will next show that there exists an interval $[T_1,T_2] \subset (T,T+\eps)$ such that
$$
\inf_{t \in [T_1,T_2]} \delta(t) > 0.
$$
Assume the contrary, i.e., that for every $t' \in (T,T+\eps)$ and for every $\theta \in(0,t'-T)$ there exists a sequence $\{t_i\}_i \subset [t'-\theta,t']$ such that $\delta(t_i) \to 0$. By repeating this with smaller and smaller $\theta$ we discover that for every $t' \in (T,T+\eps)$ there exists a sequence $\{t_i\}_i$ such that $t_i < t'$ and $t_i \to t'$ and $\delta(t_i) \to 0$ as $i \to \infty$.

Let $r=r^*(t,t')$ be the intersection of $U(\cdot,t)$ and $\omega_C$ that coincides with $r = \delta(t')$ for $t = t'$, that is,
$$
U(r^*(t,t'),t) = \omega_C(r^*(t,t')),
$$
and $r^*(t',t') = \delta(t')$. Then for any $t_0 \in (T,T+\eps)$, by the continuity of $U(\cdot,t)$ outside the origin, we can take $\eps_1, r_1>0$ small enough such that $r^*(t,t_0) > r_1>0$ for every $t \in [t_0-\eps_1,t_0]$. However, for some $t_1 \in (t_0-\eps_1,t_0)$ close to $t_0$, we have that $\delta(t_1)$ is close to $0$. Hence we may assume that $\delta(t_1) < \frac{r_1}{2}$ and we can take $\eps_2,r_2>0$ small enough so that $r^*(t,t_1) \in (r_2, \frac{r_1}{2})$ for every $t \in [t_1-\eps_2,t_1]$.

By continuing this, we know that for some $t_2\in (t_1-\eps_2,t_1)$ close to $t_1$, it holds that $\delta(t_2) \in (0,\frac{r_2}{2})$ and so $r^*(t,t_2) \in (r_3, \frac{r_2}{2})$ for every $t \in [t_2-\eps_3,t_2]$ and for some $\eps_3>0$ and $r_3 < \frac{r_2}{2}$.

We have therefore obtained a sequence $\{t_i\}_i$ such that $t_{i+1} \in (t_i-\eps_{i+1},t_i)$ and such that $U(r,t)$ intersects with $\omega_C(r)$ at $r = r^*(t,t_i) \in (r_{i+1}, \frac{r_i}{2})$ for $t \in [t_i-\eps_{i+1},t_i]$. Thus $t_i$ converges along a subsequence to a limit $t'$ such that $t' \in (t_i-\eps_{i+1},t_i)$ for every $i$ and so $U(r,t')$ intersects with $\omega_C(r)$ at $r=r^*(t',t_i)$ for every $i$. Moreover, $0 < r^*(t',t_i) < \frac{r_i}{2} < \frac{r_1}{2^i} \to 0$ as $i \to \infty$, which is in contradiction with the results we proved above, namely, the fact that $U(r,t')$ intersects $\omega_C(r)$ only finitely many times near the origin.

This shows that there exists $T_1,T_2 \in (T,T+\eps)$ and $\delta_0>0$ such that $u(x,t) < \omega_C(x)$ for $|x| < \delta_0$ and for every $t \in [T_1,T_2]$. Since $U \in C^{\alpha}_{\text{loc}}(0,R)$ for every $t \in (0,\mt)$, we know that $U$ is bounded in $[\delta_0,R_1] \times [T_1,T_2]$ for every $r_1 \in (\delta_0,R)$. The blow-up set in known to be a compact set of $B(R)$ and so $U$ is bounded in $[R_1,R] \times [T_1,T_2]$ for some $R_1<R$ which finishes the proof. \nelio

By the estimate given by Lemma \ref{pro:loglog}, we are able to use the same techniques that are used in \cite{FMP} to prove the following Lemma \ref{lem:C2reg}. The idea of the proof is to use the estimate from Lemma \ref{pro:loglog} to get that $\av u_n(\cdot,t) \av_{L^q} \le C$ and $\av e^{u_n(\cdot,t)} \av_{L^q} \le C$ for every $t \in [T_1,T_2]$ with some constant $C$ independent of $n$, where $q \in (1,\frac{N}{2})$. Similar estimates as in Lemma \ref{lemma_C_alpha} then give that $\{u_n\}_n$ is a compact subset of $H^1_0(B(R))$. This and parabolic regularity give the first part of Lemma \ref{lem:C2reg}.

To prove that $u$ is in $C^1$ with respect to $t$ we proceed as in Proposition 2.14 in \cite{FMP} and use standard parabolic estimates, Sobolev embeddings and the estimate obtained in the previous Lemma. We omit the details of the proof and refer instead to Proposition 2.12 and 2.14 in \cite{FMP}.
\begin{lem} \label{lem:C2reg}
Let $u$ be a minimal limit $L^1$-solution of (\ref{eq1})-(\ref{f(u)}) on $(0,\mt)$ that blows up at $t = T < \mt$. Then there exist $T_1,T_2 \in (T,\mt)$ such that
$$
u_n(\cdot,t) \to u(\cdot,t) \qquad \text{ as } n \to \infty \text{ in } \in H^1_0(B(R))
$$
and
$$
u_n(\cdot,t) \to u(\cdot,t) \qquad \text{ as } n \to \infty \text{ in } C^2_{\text{loc}}(B(R)\setminus \{0\}),
$$
for every $t \in (T_1,T_2)$. Here $\{u_n\}_n$ is a sequence of classical solutions defining $u$.

Moreover, $u \in C^1([T_1,T_2], L^q(B(R)))\cap C([T_1,T_2],W^{2,q}(B(R)))$, for any $q \in (1,\frac{N}{2})$.
\end{lem}

The next Lemma is proved in \cite{FMP} by assuming also that $u$ is radially nonincreasing and using the zero number diminishing property for the time derivative of $u$. Here we obtain the same result by using a technique of counting the intersections of the approximating functions $u_n$ and the singular solution
$$
\Phi^*(|x|) = \varphi^*(x) = -2\log|x| + \log(2(N-2)).
$$
\begin{lem} \label{lem:BUtimes}
Let $u$ be as in Lemma \ref{pro:loglog} and assume that $z_{[0,R]}(U_0 - \Phi^*) =M_0$. Then $u$ can blow-up at most $\frac{M_0}{2} + 1$ times with type I rate.
\end{lem}
\emph{Proof.} Assume that $u$ blows up at $t \in \{t_i\}_{i=1}^{k-1}$ incompletely with type I rate. This means in particular that $u(\cdot,t) \in L^{\infty}(B(R))$ in some interval $t \in (t_i-\eps_i,t_i)$ and hence it is a smooth function in that same interval. Therefore, the convergence $u_n(\cdot,t) \to u(\cdot,t)$ holds in the sense of $C^2(B(R))$ for $t \in (t_i-\eps_i,t_i)$.

Let $m(t_i)$ be the number of intersections of $U(\cdot,t)$ and $\Phi^*(\cdot)$ that tend to zero as $t \to t_i$. The assumptions imply that the convergence (\ref{eq:ssconv1})-(\ref{stationaryAsympt}) holds. Moreover, for every solution $\varphi(y)=\Phi(|y|)$ of (\ref{eq:stationary})-(\ref{stationaryAsympt}) one has that $z_{[0,\infty)}(\Phi-\Phi^*) \ge 2$. Consequently, $m(t_i) \ge 2$. 

By the $C^2((0,R))$ convergence of $U_n(\cdot,0) \to U(\cdot,0)$ and by the zero number diminishing property, we have that $z_{[0,R]}(U_n(\cdot,t)-\Phi^*(\cdot)) \le M_0$ for every $t>0$ and $n$ large enough. Since $U_n$ is bounded for every $t$ and because $U_n(r,t) \le U(r,t)$, we conclude that at least two intersections of $U_n(\cdot,t)$ and $\Phi^*$ vanish already for $t<t_1$ and so $U_n(\cdot,t_1)$ and $\Phi^*$ intersect at most $M_0 -2$ times, for $n$ large enough. Again using the zero number diminishing property we get that $z_{[0,R]}(U_n(\cdot,t)-\Phi^*(\cdot)) \le M_0-2$ for every $t > t_1$ and $n$ large. Moreover, since $U_n \to U$ in $C^2((0,R))$ for $t \in (t_2-\eps_2,t_2)$ we know that also $U(\cdot,t)$ and $\varphi^*$ intersect at most $M_0-2$ times in that interval.

If $u$ blows up $k$ times at $t \in \{t_i\}_{i = 1}^k$, it means that it has to blow-up incompletely at least $k-1$ times and so by a continuation of the above reasoning,
$$
z_{[0,R]}(U_n(\cdot,t_i) - \Phi^*(\cdot)) \le M_0-2i.
$$
This gives us the inequality $M_0 - 2(k-1) \ge 0$ and so $k \le \frac{M_0}{2} + 1$. \nelio

Even though we do not assume that $u$ is radially nonincreasing, the next Lemma states that if $u$ does not gain regularity immediately after the blow-up, then it is radially decreasing close to $x=0$ at least on some time interval. 
\begin{lem} \label{lem:decreasing}
Assume that $u$ is a minimal $L^1$-solution of (\ref{eq1})-(\ref{f(u)}) on $(0,\mt)$ that blows up at $t=T<\mt$ and that $\av u(\cdot,t) \av_{L^{\infty}}=\infty$ for $t \in (T,T+\eps)$. Then there exists $T_1,T_2 \in (T, T+\eps)$ and $\delta > 0$ such that
$$
U_r(r,t) < 0, \qquad \text{ for } r \in (0,\delta) \text{ and } t \in (T_1,T_2).
$$
\end{lem}

\emph{Proof.} Assume that there exists $t^* \in (T,\mt)$ such that $U_r(\cdot,t^*)$ has infinitely many zeros. Because $u_n(\cdot,t) \to u(\cdot,t)$ in $C^2(B(R)\setminus \{0\})$, we know that the number of zeros, $M(n)$, of $(U_n)_r(\cdot,t^*)$ satisfies $M(n) \to \infty$ as $n \to \infty$. Since $U_r(\cdot,t_0) \in C^{\infty}(B(R))$ has $M_0<\infty$ zeros for any $t_0<T$, we know that $U_n(\cdot,t_0)$ has $M_0$ zeros for $n$ large enough. Therefore, if $n$ is large enough, we have that $M(n) > M_0$ and there exists $t_1 \in (t_0,t^*)$ such that
\be \label{z0R_un}
z_{[0,R]}((U_n)_r(\cdot,t)) < z_{[0,R]}((U_n)_r(\cdot,t')),
\ee
for $t<t_1<t'$.

Since $u_n$ is a classical solution, we can take $\eps>0$ small and $r_0>0$ such that $(U_n)_r(r,t) < 0$ (or $>0$) for $(r,t) \in r_0 \times (t_1-\eps,t_1+\eps)$ and for $(r,t) \in (0,r_0) \times t_1-\eps$. Then, by Lemma 52.18 in \cite{QS}, we know that $(U_n)_r(r,t) < 0$ (or $>0$) for $(r,t) \in (0,r_0) \times (t_1-\eps,t_1+\eps)$. Similarly we can take $r_1 > 0$ such that $(U_n)_r(r,t) < 0$ for $(r,t) \in r_1 \times (t_1-\eps,t_1+\eps)$ and $(U_n)_r(r,t) < 0$ for $(r,t) \in (r_1,R) \times t_1-\eps$. This implies, by Proposition 52.8 in \cite{QS}, that $(U_n)(r,t) < 0$ for $(r,t) \in (r_0,R) \times (t_1-\eps,t_1+\eps)$.

Define $v = (U_n)_r$. Then
$$
v_t = v_{rr} + \frac{n-1}{r}v_r - \frac{n-1}{r^2}v + e^{U_n}v,
$$ 
and $v(r,t) < 0 $ for $(r,t) \in r_0 \times (t_1-\eps,t_1+\eps)$ and for $(r,t) \in r_1 \times (t_1-\eps,t_1+\eps)$. Hence, the standard zero number property gives
$$
z_{[r_1,r_2]}((U_n)_r(\cdot,t)) \le z_{[r_1,r_2]}((U_n)_r(\cdot,t_1-\eps)),
$$
for $t \in (t_1-\eps,t_1+\eps)$.

This contradicts the fact that, by (\ref{z0R_un}), the zero number of $(U_n)_r$ increases at $t = t_1$. Hence $U_r$ has only finitely many zeros for every $t \in (T,\mt)$ and so we can take $T_1,T_2 \in (T,T+\eps)$ and $\delta > 0$ such that the first zero of $U_r(\cdot,t)$ is at $r = r_0(t) > \delta$ for every $t \in (T_1,T_2)$. This proves the claim. \nelio

The final Lemma before concluding the proof of Theorem \ref{theorem3} improves the upper bound of Lemma \ref{pro:loglog} under the assumption that $u$ stays singular for $t>T$.
\begin{lem} \label{Lemma_log}
Under the same assumptions as in Lemma \ref{lem:decreasing}, there exist constants $T_1,T_2 \in (T,T+\eps)$ and $C, \delta>0$ such that 
$$
u(x,t) \le -2\log|x| + C,
$$
for every $t \in (T_1,T_2)$ and $x \in B(\delta)$.
\end{lem}
\emph{Proof.} Because of the Lemmata \ref{lem:C2reg} and \ref{lem:decreasing}, we can use the same method as in Lemma 2.16 in \cite{FMP} to prove this Lemma so let us only give the outline of the proof.

Let $\mu$ be the eigenvalue and $\eta$ the eigenfunction of the radial Laplacian in $B(1)$ under the Neumann boundary condition, that is,
$$
\left\{
\begin{array}{ll}
\eta''(r) + \frac{N-1}{r}\eta'(r) = - \mu\eta(r), & r \in (0,1), \\
\eta(0) = 1, \, \eta'(0) = \eta'(1) = 0,
\end{array}
\right.
$$
and $\eta'(r) < 0$, for $r \in (0,1)$.

Define $\psi_{\rho}(x) = \eta(|x|/\rho) - \eta(1)$, which implies
$$
\left\{
\begin{array}{ll}
\Delta \psi_{\rho} = -\frac{\mu}{\rho^2}(\psi_{\rho} + \eta(1)), & x \in B(\rho), \\
\psi_{\rho}> 0, & x \in B(\rho), \\
\psi_{\rho} = \frac{\partial \psi_{\rho}}{\partial n} = 0, & x \in \partial B(\rho).
\end{array}
\right.
$$
Then let
$$
h_{\rho}(t) = \left(\int_{B(\rho)} \psi_{\rho}(x) \ud x \right)^{-1} \int_{B(\rho)} \psi_{\rho}(x) u(x,t) \ud x,
$$
and use Lemma \ref{lem:C2reg} to differentiate $h_{\rho}$ and estimate, as in \cite{FMP}, to obtain
$$
h_{\rho}'(t) \ge e^{h_{\rho}(t)} - \frac{\mu}{\rho^2} h_{\rho}(t) + \frac{\mu}{\rho^2} \frac{1}{|B(\rho)|} \int_{B(\rho)} u(x,t) \ud x.
$$

By the previous Lemma, we have $T_1,T_2 \in (T,T+\eps)$ and $\delta > 0$ such that $U_r(r,t) < 0$ for $r \in (0,\delta)$ and $t \in (T_1,T_2)$. By the zero number diminishing property and Proposition \ref{pro:Va}, we may also assume that 
$$
U(r,t) \ge \Phi^*(r) = -2\log(r) + \log(2(N-2)),
$$
for $r \in (0,\delta)$ and $t \in (T_1,T_2)$. Therefore, for $\rho < \delta$, one has
$$
h_{\rho}'(t) \ge e^{h_{\rho}(t)} - \frac{\mu}{\rho^2} \Big(h_{\rho}(t) +2\log(\rho) - \log(2(N-2))\Big),
$$
and by defining $k(t) = h_{\rho}(t) +2\log(\rho)$, we get
$$
k'(t) \ge \frac{1}{\rho^2}\Big(e^{k(t)} - \mu k(t) + \mu \log(2(N-2))\Big).
$$
Integrating from $t \in (T_1,T_2)$ to $T_2$ implies
$$
\int_{k(t)}^{k(T_2)} \frac{\ud k}{e^k - \mu k + \mu \log(2(N-2))} \ge \frac{T_2-t}{\rho^2} > \frac{T_2-T_1}{2R^2},
$$
for every $\rho < \delta$ and $t\in (T_1,\frac{T_1+T_2}{2})$. Therefore, either $k(t) < k^*$, where $k^*$ is the largest root of $e^k - \mu k + \mu \log(2(N-2))$, or 
$$
\int_{k(t)}^{\infty} \frac{\ud k}{e^k - \mu k + \mu \log(2(N-2))} \ge \frac{T_2-T_1}{2R^2}.
$$
Defining $a>0$ to be such that
$$
\int_a^{\infty} \frac{\ud k}{e^k - \mu k + \mu \log(2(N-2))} = \frac{T_2-T_1}{2R^2},
$$
one thus have $k(t) < a$. Since $U_r(r,t)<0$ for $r < \delta$ and $t \in (T_1,T_2)$, we have $U(r,t) < h_r(t)$ for $r < \delta$ and so the claim is proved. \nelio

Theorem \ref{theorem3} can now be proved by using the above results and the technique from \cite{FMP}.
\begin{pro} \label{pro:reg}
Let the assumptions of Theorem \ref{theorem3} be satisfied. Then $u$ becomes regular immediately after the blow-up, i.e., there exists $\eps>0$ such that $u \in C^{\infty}(B(R))$ for $t \in (T,T+\eps)$. 
\end{pro}
\emph{Proof.} Assume that the claim is false. Then there exists a sequence $\{t_i\}_i$ tending to $T$ from above such that $\av u(\cdot,t_i) \av_{L^{\infty}} = \infty$ for every $i$. If there exists another sequence $\{\tau_i\}_i$ tending to $T$ from above such that $\av u(\cdot,\tau_i) \av_{L^{\infty}} < \infty$, then $u$ is regular in some interval $t \in (\tau_i,\tau_i + \eps_i)$ and so $u$ would blow-up infinitely many times after $t = T$. This contradicts Lemma \ref{lem:BUtimes} and therefore $\av u(\cdot,t) \av_{L^{\infty}} = \infty$ for $t \in (T,T+\eps)$ and for some $\eps>0$.

Define, for some $\theta \in (T_1,T_2)$,
$$
w_{\theta}(y,s) = -\log(\theta-t) + u(\sqrt{\theta-t}\,y,t),
$$
where $s = -\log(\theta-t)$. Then by Lemma \ref{Lemma_log} and Proposition \ref{pro:Va}, and since $u(\cdot,t)\in C((\delta,R))$, we have constants $C_1$, $\delta > 0$ and $T_1,T_2 \in (T,T+\eps)$ such that
\be \label{u_strict_est}
-2\log|x| + \log(2(N-2)) \le u(x,t) \le -2\log|x| + C_1,
\ee
for $|x| \in (0,\delta)$ and $t \in (T_1,T_2)$. This gives
$$
-2\log|y|- C_1' \le w_{\theta}(y,s) \le -2\log|y| + C_2',
$$
for some $C_1'$ and $C_2'$ and $s$ large enough. Then an energy argument implies that $w_{\theta}(\cdot,s)$ converges, along a subsequence $\{s_i\}_i$ tending to infinity, to a limit $\varphi$ that satisfies
$$
\Delta \varphi - \frac{y}{2}\nabla \varphi + e^{\varphi}-1=0,
$$
with asymptotics (\ref{stationaryAsympt}).

By (\ref{u_strict_est}), it holds that $\varphi(y) \ge \varphi^*(y) = -2\log|y| + \log(2(N-2))$. Since there are no such $\varphi$ that would lie above $\varphi^*$, $w_{\theta}(\cdot,s)$ has to converge to $\varphi^*$.

However, $w_{\theta}$ can not converge to $\varphi^*$, as can be demonstrated by an intersection number argument. This finishes the proof. We refer to \cite{FMP} for the details. \nelio

\section{Blow-up rate for subcritical dimensions} \label{section_BUrate}
The aim of this section is to prove Theorem \ref{theo:n=2}. The proof follows the lines of \cite{FP} and \cite{MM1}, where type I blow-up is proved for supercritical and critical cases.

In \cite{FP} it is proved that blow-up is of type I in the case of the exponential nonlinearity for $N \in (2,10)$, provided that the solution is radially symmetric and attains its maximum at the origin. The following Lemma, proved in \cite{FP}, is independent of the value of the dimension $N$.

Let us again use the radial notation, i.e., $U(|x|,t) = u(x,t)$.
\begin{lem} \label{lemma_conv}
Let $u$ be a radially symmetric solution of (\ref{eq1})-(\ref{f(u)}) that blows up at $t = T$ and satisfies $u(0,t) = \max_x u(x,t)$ for every $t \in (0,T)$. Let
$$
w_i(\rho,\tau) = -u(0,t_i) + U(e^{-u(0,t_i)/2}\rho, t_i + e^{-u(0,t_i)}\tau),
$$
for $(\rho,\tau) \in (0,u^{u(0,t_i)/2}R) \times (-e^{u(0,t_i)}t_i,e^{u(0,t_i)}(T-t_i))$. 

If blow-up is of type II, then there exists a sequence $\{t_i\}_i$ tending to $T$ such that
$$
w_i(\rho,\tau) \to \psi(\rho),
$$
uniformly for $(\rho,\tau)$ in compact sets of $(0,\infty)\times (-\infty,\frac{1}{4})$, where 
\be \label{steadystates}
\left\{
\begin{array}{ll}
\psi'' + \frac{N-1}{\rho}\psi' + e^{\psi} = 0,& \text{ for } \rho>0, \\
\psi(0) = 0, \, \psi'(0) = 0.
\end{array}
\right.
\ee
\end{lem}

Theorem \ref{theo:n=2} can now be proved by a technique from \cite{MM1}, where type I blow-up was proved for the power nonlinearity when $p  = p_S$ and $u$ is radially symmetric.

\emph{Proof of Theorem \ref{theo:n=2}.} By results in \cite{T}, we have that the only solutions of (\ref{steadystates}) are
$$
\left\{
\begin{array}{ll}
\psi(\rho) = -2\log\left( \cosh\left(\frac{r}{\sqrt{2}}\right) \right),& \text{for } N=1,\\
\psi(\rho) = -2\log\left(1 + \frac{r^2}{8} \right), & \text{for } N = 2.
\end{array}
\right.
$$
By the scaling invariance of the equation, we know that $\psi_a(\rho) = a + \psi(e^{a/2}\rho)$ is also a solution of (\ref{steadystates}) with $\psi_a(0) = a$, for every $a \in \mbr$. Furthermore, $\psi_a$ crosses $\psi_b$ at least once for every $a \ne b$.

Assume that blow-up is of type II and let $\{t_i\}_i$ and $w_i$ be as in the previous Lemma. Then, by the zero number diminishing property, for every $t_0 \in (0,T)$, there exists $a_0(t_0)$ such that
$$
z_{[0,\frac{R}{2}]}(U(\cdot,t)-\psi_a(\cdot)) \le 1,
$$
for every $t \ge t_0$.

Since $U(0,t) \to \infty$, as $t \to T$, one notices that, for every $a>0$, there exists $t(a) \in (0,T)$ such that $z_{[0,\frac{R}{2}]}(U(\cdot,t)-\psi_a(\cdot)) = 0$ for $t>t(a)$. Moreover, $t(a)$ can be taken to be the time moment for which $U(0,t(a)) = \psi_a(0) = a$. This implies that $U(r,t) > \psi_a(r)$ for $r \in [0,\frac{R}{2}]$ and $t > t(a)$.

By the previous Lemma, one has
$$
-u(0,t_i) + U(e^{-u(0,t_i)/2} r,t_i) \to \psi_0(r),
$$
as $i \to \infty$. Hence,
$$
U(0,t_i) > u(0,t_i)-1 = \psi_{u(0,t_i)-1}(0),
$$
for $i$ large enough, which gives that $t_i > t(u(0,t_i)-1)$ for $i$ large enough. This allows us to conclude that
$$
U(r,t_i) > \psi_{u(0,t_i)-1}(r),
$$
for every $r \in [0,\frac{R}{2}]$ and $i$ large enough. By scaling, this translates into
$$
-u(0,t_i) + U(e^{-u(0,t_i)/2}r,t_i) > -u(0,t_i) + \psi_{u(0,t_i)-1}(e^{-u(0,t_i)/2}r) = \psi_{-1}(r),
$$
for every $r \in [0,R e^{u(0,t_i)/2}/2]$.

By taking the limit $i \to \infty$, the previous Lemma implies 
$$
\psi_0(r) > \psi_{-1}(r),
$$
for every $r \ge 0$. This contradicts the results in \cite{T}, where it is proved that $\psi_0$ and $\psi_{-1}$ intersect at least once. The claim thus follows. \nelio


\end{document}